\documentclass[reqno, 12pt]{amsart}
\makeatletter
\let\origsection=\section \def\section{\@ifstar{\origsection*}{\mysection}} 
\def\mysection{\@startsection{section}{1}\z@{.7\linespacing\@plus\linespacing}{.5\linespacing}{\normalfont\scshape\centering\S\hspace{1pt}}}
\makeatother

\usepackage[T1]{fontenc}
\usepackage{lmodern}
\usepackage[babel]{microtype}
\usepackage[english]{babel}

\usepackage{geometry}
\usepackage{mathdots}
\usepackage[
sorting=nty,
backend=biber,
natbib=true,
style=numeric,   
isbn=false,
maxbibnames=9,
backref
]{biblatex}
\usepackage{csquotes}
\DeclareFieldFormat{title}{\mkbibquote{#1}}
\DeclareFieldFormat[misc]{date}{Preprint (#1)}

\renewbibmacro*{doi+eprint+url}{
	\iftoggle{bbx:url}
	{\iffieldundef{doi}{\usebibmacro{url+urldate}}{}}{}
	\newunit\newblock
	\iftoggle{bbx:eprint}
	{\usebibmacro{eprint}}
	{}
	\newunit\newblock
	\iftoggle{bbx:doi}
	{\printfield{doi}}
	{}
}

\renewbibmacro{in:}{}
\DeclareFieldFormat{pages}{#1}

\DefineBibliographyStrings{english}{  backrefpage={$\uparrow$}, 
backrefpages={$\uparrow$}}
\usepackage{csquotes}

\usepackage{mathrsfs}
\usepackage[colorlinks=true, linkcolor=blue, filecolor=magenta, urlcolor=cyan]{hyperref}
\usepackage{xcolor}
\hypersetup{
    colorlinks,
    linkcolor={red!80!black},
    citecolor={red!60!black},
    urlcolor={blue!60!black}
 }

\usepackage{amsmath,amssymb,amsthm}
\usepackage{mathrsfs}
\usepackage{mathabx}\changenotsign
\usepackage{dsfont}

\numberwithin{equation}{section}
\numberwithin{figure}{section}

\usepackage[shortlabels]{enumitem}

\def\caselabel{\textit{{{\hskip 0.7em Case~}}\upshape(\!{\itshape \arabic*\,}\!):}}

\newcommand{\customlabel}[2]{#2\def\@currentlabel{#2}\label{#1}}
\def\rmlabel{\upshape({\itshape \roman*\,})}

\usepackage{tikz}
\usepackage{scalerel}

\usetikzlibrary{calc,decorations.pathmorphing}
\pgfdeclarelayer{background}
\pgfdeclarelayer{foreground}
\pgfdeclarelayer{front}
\pgfsetlayers{background,main,foreground,front}

\usepackage{subcaption}

\newcommand{\qedge}[7]{

	\ifx\relax#4\relax
		\def\qoffs{0pt}
	\else
		\def\qoffs{#4}
	\fi

	\def\qhedge{
		($#1+#3!\qoffs!-90:#2-#3$) --
		($#2+#1!\qoffs!-90:#3-#1$) --
		($#3+#2!\qoffs!-90:#1-#2$) -- cycle}

	\coordinate (12) at ($#1!\qoffs!90:#2$);
	\coordinate (13) at ($#1!\qoffs!-90:#3$);
	\coordinate (23) at ($#2!\qoffs!90:#3$);
	\coordinate (21) at ($#2!\qoffs!-90:#1$);
	\coordinate (31) at ($#3!\qoffs!90:#1$);
	\coordinate (32) at ($#3!\qoffs!-90:#2$);
	
	\def\nqhedge{
		(13) let \p1=($(13)-#1$), \p2=($(12)-#1$) in
			arc[start angle={atan2(\y1,\x1)}, delta angle={atan2(\y2,\x2)-atan2(\y1,\x1)-360*(atan2(\y2,\x2)-atan2(\y1,\x1)>0)}, x radius=\qoffs, y radius=\qoffs] --
		(21) let \p1=($(21)-#2$), \p2=($(23)-#2$) in
			arc[start angle={atan2(\y1,\x1)}, delta angle={atan2(\y2,\x2)-atan2(\y1,\x1)-360*(atan2(\y2,\x2)-atan2(\y1,\x1)>0)}, x radius=\qoffs, y radius=\qoffs] --
		(32) let \p1=($(32)-#3$), \p2=($(31)-#3$) in
			arc[start angle={atan2(\y1,\x1)}, delta angle={atan2(\y2,\x2)-atan2(\y1,\x1)-360*(atan2(\y2,\x2)-atan2(\y1,\x1)>0)}, x radius=\qoffs, y radius=\qoffs] --
		cycle}

		\ifx\relax#5\relax
		\def\qlwidth{1pt}
	\else
		\def\qlwidth{#5}
	\fi
	
		\ifx\relax#7\relax
		\fill \nqhedge;
	\else
		\fill[#7]\nqhedge;
	\fi

		\ifx\relax#6\relax
		\draw[line width=\qlwidth,rounded corners=\qoffs]\nqhedge;
	\else
		\draw[line width=\qlwidth,#6]\nqhedge;
	\fi
}

\newcommand{\redge}[8]{

		\ifx\relax#5\relax
		\def\qoffs{0pt}
	\else
		\def\qoffs{#5}
	\fi

				\def\rhedge{
		($#1+#4!\qoffs!-90:#2-#4$) -- 
		($#2+#1!\qoffs!-90:#3-#1$) -- 
		($#3+#2!\qoffs!-90:#4-#2$) -- 
		($#4+#3!\qoffs!-90:#1-#3$) -- cycle}

	\coordinate (12) at ($#1!\qoffs!90:#2$);
	\coordinate (14) at ($#1!\qoffs!-90:#4$);
	\coordinate (23) at ($#2!\qoffs!90:#3$);
	\coordinate (21) at ($#2!\qoffs!-90:#1$);
	\coordinate (34) at ($#3!\qoffs!90:#4$);
	\coordinate (32) at ($#3!\qoffs!-90:#2$);
	\coordinate (41) at ($#4!\qoffs!90:#1$);
	\coordinate (43) at ($#4!\qoffs!-90:#3$);
	
	\def\nrhedge{
		(14) let \p1=($(14)-#1$), \p2=($(12)-#1$) in 
			arc[start angle={atan2(\y1,\x1)}, delta angle={atan2(\y2,\x2)-atan2(\y1,\x1)-360*(atan2(\y2,\x2)-atan2(\y1,\x1)>0)}, x radius=\qoffs, y radius=\qoffs] --
		(21) let \p1=($(21)-#2$), \p2=($(23)-#2$) in 
			arc[start angle={atan2(\y1,\x1)}, delta angle={atan2(\y2,\x2)-atan2(\y1,\x1)-360*(atan2(\y2,\x2)-atan2(\y1,\x1)>0)}, x radius=\qoffs, y radius=\qoffs] --
		(32) let \p1=($(32)-#3$), \p2=($(34)-#3$) in 
			arc[start angle={atan2(\y1,\x1)}, delta angle={atan2(\y2,\x2)-atan2(\y1,\x1)-360*(atan2(\y2,\x2)-atan2(\y1,\x1)>0)}, x radius=\qoffs, y radius=\qoffs] --
		(43) let \p1=($(43)-#4$), \p2=($(41)-#4$) in 
			arc[start angle={atan2(\y1,\x1)}, delta angle={atan2(\y2,\x2)-atan2(\y1,\x1)-360*(atan2(\y2,\x2)-atan2(\y1,\x1)>0)}, x radius=\qoffs, y radius=\qoffs] --
		cycle}

		\ifx\relax#6\relax
		\def\rlwidth{1pt}
	\else
		\def\rlwidth{#6}
	\fi
	
		\ifx\relax#8\relax
		\fill \nrhedge;
	\else
		\fill[#8]\nrhedge;
	\fi

		\ifx\relax#7\relax
		\draw[line width=\rlwidth,rounded corners=\qoffs]\nrhedge;
	\else
		\draw[line width=\rlwidth,#7]\nrhedge;
	\fi
	}

\newsavebox\vdegbox
\savebox\vdegbox{\tikz{
		\draw[black,fill=black] (90:1) circle (.35);
		\draw[black,line width=0.10cm] (210:1) circle (.30);
		\draw[black,line width=0.10cm] (330:1) circle (.30);
		\draw[opacity=0] (0:1.2) circle (0.1);
	}}

\newsavebox\vvbox
\savebox\vvbox{\tikz{
		\draw[black,line width=0.10cm] (90:1) circle (.30);
		\draw[black,fill=black] (210:1) circle (.35);
		\draw[black,fill=black] (330:1) circle (.35);
		\draw[opacity=0] (0:1.2) circle (0.1);
	}}

\newsavebox\pdegbox
\savebox\pdegbox{\tikz{
		\draw[black,line width=0.10cm] (90:1) circle (.30);
		\draw[black,fill=black] (210:1) circle (.35);
		\draw[black,fill=black] (330:1) circle (.35);
		\draw[black,line width=0.28cm ] (210:1) -- (330:1);
		\draw[opacity=0] (0:1.2) circle (0.1);
	}}

\newsavebox\vvvbox
\savebox\vvvbox{\tikz{
		\draw[black,fill=black] (90:1) circle (.35);
		\draw[black,fill=black] (210:1) circle (.35);
		\draw[black,fill=black] (330:1) circle (.35);
		\draw[opacity=0] (0:1.2) circle (0.1);
	}}

\newsavebox\evbox
\savebox\evbox{\tikz{
		\draw[black,fill=black] (90:1) circle (.35);
		\draw[black,fill=black] (210:1) circle (.35);
		\draw[black,fill=black] (330:1) circle (.35);
		\draw[black,line width=0.28cm ] (210:1) -- (330:1);
		\draw[opacity=0] (0:1.2) circle (0.1);
	}}

\newsavebox\eebox
\savebox\eebox{\tikz{
		\draw[black,fill=black] (90:1) circle (.35);
		\draw[black,fill=black] (210:1) circle (.35);
		\draw[black,fill=black] (330:1) circle (.35);
		\draw[black,line width=0.28cm ] (90:1) -- (330:1);
		\draw[black,line width=0.28cm ] (90:1) -- (210:1);
		\draw[opacity=0] (0:1.2) circle (0.1);
	}}

\newsavebox\eeebox
\savebox\eeebox{\tikz{
		\draw[black,fill=black] (90:1) circle (.35);
		\draw[black,fill=black] (210:1) circle (.35);
		\draw[black,fill=black] (330:1) circle (.35);
		\draw[black,line width=0.28cm ] (90:1) -- (330:1);
		\draw[black,line width=0.28cm ] (90:1) -- (210:1);
		\draw[black,line width=0.28cm ] (210:1) -- (330:1);
		\draw[opacity=0] (0:1.2) circle (0.1);
	}}

\theoremstyle{plain}
\newtheorem{thm}{Theorem}[section]
\newtheorem{fact}[thm]{Fact}
\newtheorem{prop}[thm]{Proposition}
\newtheorem{clm}{Claim}
\newtheorem{prob}[thm]{Problem}
\newtheorem{quest}[thm]{Question}
\newtheorem{cor}[thm]{Corollary}
\newtheorem{lemma}[thm]{Lemma}

\theoremstyle{definition}

\newtheorem{dfn}[thm]{Definition}

\newenvironment{claimproof}[1]{\par\noindent\textit{Proof of the claim:}\space#1}{\leavevmode\unskip\penalty9999 \hbox{}\nobreak\hfill\quad\hbox{$\blacksquare$}\smallskip}

\def\back#1{\ThisStyle{%
  \mathord{\vbox{\offinterlineskip\ialign{
    \hfil##\hfil\cr
    $\SavedStyle{}_{\smash{\,\leftharpoonup}}$\cr
    \noalign{\kern-0.5\scriptspace}
    $\SavedStyle#1$\cr}}}}}
\def\forw#1{\ThisStyle{%
  \mathord{\vbox{\offinterlineskip\ialign{
    \hfil##\hfil\cr
    $\SavedStyle{}_{\smash{\,\rightharpoonup}}$\cr
    \noalign{\kern-0.5\scriptspace}
    $\SavedStyle#1$\cr}}}}}

\let\eps=\varepsilon
\let\theta=\vartheta
\let\rho=\varrho
\let\phi=\varphi

\DeclareMathOperator{\ex}{ex}

\let\polishlcross=\l
\def\l{\ifmmode\ell\else\polishlcross\fi}

\makeatletter
\def\moverlay{\mathpalette\mov@rlay}
\def\mov@rlay#1#2{\leavevmode\vtop{%
   \baselineskip\z@skip \lineskiplimit-\maxdimen
   \ialign{\hfil$\m@th#1##$\hfil\cr#2\crcr}}}
\newcommand{\charfusion}[3][\mathord]{
    #1{\ifx#1\mathop\vphantom{#2}\fi
        \mathpalette\mov@rlay{#2\cr#3}
      }
    \ifx#1\mathop\expandafter\displaylimits\fi}
\makeatother

\newcommand{\dcup}{\charfusion[\mathbin]{\cup}{\cdot}}

\def\qand{\quad\text{and}\quad}
\def\qqand{\qquad\text{and}\qquad}
\def\nquad{\!\!\!\!}
\def\nqquad{\!\!\!\!\!\!\!\!}

\newcommand{\vrhup}[1]{\scaleobj{0.6}{\scalerel*{\rightharpoonup}{#1}}}
\newcommand{\nrhup}{\mathord{\scaleobj{0.6}{\scalerel*{\rightharpoonup}{x}}}}

\def\vseq#1{\ThisStyle{%
  \mathord{\vbox{\offinterlineskip\ialign{%
    \hfil##\hfil\cr
    $\SavedStyle{}_{\smash{\vrhup#1}}$\cr
    \noalign{\kern-0.7\scriptspace}
    $\SavedStyle#1$\cr}}}}}

\def\seq#1{\ThisStyle{%
  \mathord{\vbox{\offinterlineskip\ialign{%
    \hfil##\hfil\cr
    $\SavedStyle{}_{\smash{\nrhup}}$\cr
    \noalign{\kern-0.5\scriptspace}
    $\SavedStyle#1$\cr}}}}}

\let\setminus=\smallsetminus
\let\emptyset=\varnothing

\let\to=\lra

\usepackage{xspace}

\makeatletter
\newcommand{\pushright}[1]{\ifmeasuring@#1\else\omit\hfill$\displaystyle#1$\fi\ignorespaces}
\newcommand{\pushleft}[1]{\ifmeasuring@#1\else\omit$\displaystyle#1$\hfill\fi\ignorespaces}
\makeatother

\title[Arbitrarily small codegree Tur\'an density]{Hypergraphs with arbitrarily small codegree Tur\'an density}

\author[S. Piga]{Sim\'on Piga}
\address[S. Piga]{School of Mathematics, University of Birmingham, Edgbaston, Birmingham, B15 2TT, UK}
\email{s.piga@bham.ac.uk}
\author[B. Sch\"ulke]{Bjarne Sch\"ulke}
\address[B. Sch\"ulke]{Mathematics Department, California Institute of Technology, USA}
\email{schuelke@caltech.edu}
\thanks{
The research leading to these results was partially supported by EPSRC, grant no. EP/V002279/1 (S.~Piga).
There are no additional data beyond that contained within the main manuscript.}

\thanks{}

\subjclass[]{}

\begin{document}

\begin{abstract}
    Let $k\ge 3$.
    Given a $k$-uniform hypergraph $H$, the minimum codegree~$\delta(H)$ is the largest $d\in\mathds N$ such that every $(k-1)$-set of $V(H)$ is contained in at least $d$ edges.
    Given a $k$-uniform hypergraph $F$, the codegree  Tur\'an density $\gamma(F)$ of~$F$ is the smallest~$\gamma \in [0,1]$ such that every $k$-uniform hypergraph on $n$ vertices with $\delta(H)\geq (\gamma + o(1))n$ contains a copy of $F$.
    Similarly as other variants of the hypergraph Tur\'an problem, determining the codegree Tur\'an density of a hypergraph is in general notoriously difficult and only few results are known.
    
    In this work, we show that for every~$\eps>0$, there is a~$k$-uniform hypergraph~$F$ with~$0<\gamma(F)<\eps$.
    This is in contrast to the classical Tur\'an density, which cannot take any value in the interval $(0,k!/k^k)$ due to a fundamental result by Erd\H os.
\end{abstract} 

\maketitle

\section{Introduction}

    A $k$-uniform hypergraph (or~$k$-graph) $H$ consists of a vertex set~$V(H)$ together with a set of edges~$E(H)\subseteq V(H)^{(k)}=\{S\subseteq V(H):\vert S\vert =k\}$.
	Given a~$k$-graph~$F$ and~$n\in\mathds{N}$, the Tur\'an number of~$n$ and~$F$, $\ex(n,F)$, is the maximum number of edges an~$n$-vertex~$k$-graph can have without containing a copy of~$F$. 
    Since the main interest lies in the asymptotics, the \emph{Tur\'an density}~$\pi(F)$ of a~$k$-graph~$F$ is defined as
    $$\pi(F) = \displaystyle\lim_{n \to \infty} \frac{\ex(n,F)}{\binom{n}{k}}\,.$$
	Determining the value of~$\pi(F)$ for~$k$-graphs (with~$k\geq 3$) is one of the central open problems in combinatorics. 
	In particular, the problem of determining the Tur\'an density of the complete~$3$-graph on four vertices, i.e.,~$\pi(K_4^{(3)})$, was asked by Tur\'an in 1941~\cite{T:41} and Erd\H{o}s~\cite{E:77} offered 1000\$ for its resolution.
	Despite receiving a lot of attention (see for instance the survey by Keevash~\cite{K:11}), this problem, and even the seemingly simpler problem of determining~$\pi(K_4^{(3)-})$, where~$K_4^{(3)-}$ is the~$K_4^{(3)}$ minus one edge, remain open.
	
	Several variations of this type of problem have been considered, see for instance~\cites{BCL:21,Erdossos,Christiansurvey} and the references therein.
	The variant that we are concerned with here asks how large the \emph{minimum codegree} of an~$F$-free~$k$-graph can be.
	Given a~$k$-graph~$H=(V,E)$ and~$S\subseteq V$, the degree~$d(S)$ of~$S$ (in~$H$) is the number of edges containing~$S$, i.e., $d(S)=\vert\{e\in E:S\subseteq e\}\vert$.
	The \emph{minimum codegree} of~$H$ is defined as~$\delta(H)=\min_{x\in V^{(k-1)}}d(x)$.
	
	Given a~$k$-graph~$F$ and~$n\in\mathds{N}$, Mubayi and Zhao \cite{MZ:07} introduced the \emph{codegree Tur\'an number}~$\ex_{\text{co}}(n,F)$ of~$n$ and~$F$ as the maximum~$d$ such that there is an~$F$-free $k$-graph~$H$ on~$n$ vertices with~$\delta(H)\geq d$. 
	Moreover, they defined the \emph{codegree Tur\'an density}~$\gamma(F)$ of $F$ as
	$$\gamma(F) = \lim_{n\to\infty} \frac{ex_{\text{co}}(n,F)}{n}\,$$
    and proved that this limit always exists. 
    It is not hard to see that~$\gamma(F) \leq \pi(F)\,.$
    The codegree Tur\'an density of a family~$\mathcal{F}$ of~$k$-graphs is defined analogously.

    Similarly as for the Tur\'an density, determining the exact codegree Tur\'an density of a given hypergraph can be very difficult and so it is only known for very few hypergraphs (see the table in~\cite{BCL:21}).

    In this work, we show that there are~$k$-graphs with arbitrarily small but strictly positive codegree Tur\'an densities.

    \begin{thm}\label{thm:main}
        For every~$\xi>0$ and~$k\geq 3$, there is a~$k$-graph~$F$ with~$0<\gamma(F)<\xi$\,.
    \end{thm}
    
    Note that this is in stark contrast to the Tur\'an density and the uniform Tur\'an density, another variant of the Tur\'an density that was introduced by Erd\H{o}s and S\'os~\cite{Erdossos}.
    Regarding the former, a classical result by Erd\H{o}s~\cite{E:64} states that for no~$k$-graph the Tur\'an density is in the interval~$(0,k!/k^k)$.
    Regarding the latter Reiher, R\"odl, and Schacht~\cite{RRS:18} proved that for no~$3$-graph the uniform Tur\'an density is in~$(0,1/27)$.
    Mubayi and Zhao~\cite{MZ:07} defined $$\Gamma^{(k)}:=\{\gamma(F) \colon F \text{ is a~$k$-graph}\}\subseteq[0,1]$$ and $$\widetilde{\Gamma}^{(k)}:=\{\gamma(\mathcal{F}) \colon \mathcal{F} \text{ is a family of~$k$-graphs}\}\subseteq[0,1]\,.$$
    We remark that~$\Gamma^{(k)}\subseteq\widetilde{\Gamma}^{(k)}$ and that similar sets have been studied for the classical Tur\'an density (see, for instance,~\cite{BT:11,FR:84,P:14,S:23}).
    Mubayi and Zhao~\cite{MZ:07} showed that~$\widetilde{\Gamma}^{(k)}$ is dense in~$[0,1]$ and asked if this is also true for~$\Gamma^{(k)}$.
    Their proof for~$\widetilde{\Gamma}^{(k)}$ is based on showing that zero is an accumulation point of~$\widetilde{\Gamma}^{(k)}$.
    Theorem~\ref{thm:main} implies the same for $\Gamma^{(k)}$.
    
    \begin{cor}
        Zero is an accumulation point of~$\Gamma^{(k)}$.
    \end{cor}

Given a~$k$-graph~$H=(V,E)$ and a subset of vertices~$A=\{v_1,\dots, v_s\} \subseteq V$, we omit parentheses and commas and simply write~$A=v_1\cdots v_s$.
For the proof of Theorem~\ref{thm:main}, we consider the following hypergraphs.

\begin{dfn}\label{dfn:Zell}
    For integers~$\ell \geq k\geq 2$, we define the \emph{$k$-uniform zycle of length $\ell$} as the~$k$-graph~$Z_\ell^{(k)}$ given by 
    \begin{align*}
        V(Z_\ell^{(k)})=&\{v_i^j\colon i\in [\ell], j\in [k-1] \} \text{, and}\\
        E(Z_\ell^{(k)})=&\{v_i^1v_i^2\cdots v_i^{k-1} v_{i+1}^j \colon i\in [\ell], j\in[k-1]\}\,,
    \end{align*}    
    where the sum of indices is taken modulo~$\ell$.
\end{dfn}
Observe that~$Z_\ell^{(k)}$ has $(k-1)\ell$ vertices and $(k-1)\ell$ edges.
Moreover, $Z_\ell^{(2)} = C_\ell$. 
When~$k\in \mathds N$ is clear from the context, we omit it in the notation.

\begin{figure}
     \centering
     \begin{subfigure}[b]{0.45\textwidth}
     \centering
         \begin{tikzpicture}[scale=0.7]
    \foreach \a in {1,...,6}
    {
        \qedge{({\a*60}:2)}{({\a*60}:3)}{({(\a-1)*60}:2)}{4.5pt}{1.5pt}{red!70!black}{red!70!black,opacity=0.2};
        \qedge{({\a*60}:2)}{({\a*60}:3)}{({(\a-1)*60}:3)}{4.5pt}{1.5pt}{red!70!black}{red!70!black,opacity=0.2};
        \fill  ({\a*60}:2)  circle (2pt);
        \fill  ({\a*60}:3)  circle (2pt);
    }       
\end{tikzpicture}
        \caption{Copy of~$Z_6^{(3)}$}
\end{subfigure}
\hfill
    \begin{subfigure}[b]{0.45\textwidth}
    \centering
         \begin{tikzpicture}[scale=0.55]
        \foreach \a in {1,...,8}
        {
        \redge{({\a*45}:2)}{({\a*45}:3)}{({\a*45}:4)}{({(\a-1)*45}:2)}{4.5pt}{1.5pt}{yellow!80!black}{yellow,opacity=0.2};
        \redge{({\a*45}:2)}{({\a*45}:3)}{({\a*45}:4)}{({(\a-1)*45}:3)}{4.5pt}{1.5pt}{yellow!80!black}{yellow,opacity=0.2};
        \redge{({\a*45}:2)}{({\a*45}:3)}{({\a*45}:4)}{({(\a-1)*45}:4)}{4.5pt}{1.5pt}{yellow!80!black}{yellow,opacity=0.2};
        \fill  ({\a*45}:2)  circle (2pt);
        \fill  ({\a*45}:3)  circle (2pt);
        \fill  ({\a*45}:4)  circle (2pt);
        }       
        \end{tikzpicture}
        \caption{Copy of~$Z_8^{(4)}$}
    \end{subfigure}
\end{figure}

The following bounds on the codegree Tur\'an density of zycles imply Theorem~\ref{thm:main}.

\begin{thm}\label{thm:main2}
    Let~$k\geq 3$. For every~$d\in (0,1]$, there is an~$\ell\in \mathds N$ such that
    $$\frac{1}{2(k-1)^{\ell}}\leq \gamma (Z_\ell) \leq d\,.$$
\end{thm}

In fact we show that~$\gamma (Z_\ell)>0$ for every~$\ell\geq 3$ (see Lemma~\ref{lem:lower}).

Finally, we prove that any proper subgraph of~$Z^{(3)}_{\ell}$ has codegree Tur\'an density zero.
Let~$Z^{(3)-}_{\ell}$ be the~$3$-graph obtained from~$Z^{(3)}_{\ell}$ by deleting one edge.

\begin{thm}\label{thm:min}
Let~$\ell\geq 3$. Then
    $\gamma(Z_\ell^{(3)-})=0\,.$
\end{thm}

To prove Theorem~\ref{thm:min}, we generalise a method developed by the authors together with Sales in~\cite{PSS}.
    
\section{Proof of Theorem~\ref{thm:main2}}

Given a~$k$-graph~$H=(V,E)$, we define the neighbourhood of~$x\in V^{(k-1)}$ as $$N(x)=\{v\in V:x\cup\{v\}\in E\}\,.$$
Given a~$(k-1)$-subset of vertices~$e\in V^{(k-1)}$, we define the \emph{back neighbourhood of $e$} and the \emph{back degree of~$e$}, respectively, by
$$\back{N}(e) = \{f\in V^{(k-1)}\colon f\cup \{v\} \in E \text{ for every }v\in e\} \qand \back{d}(e)=\left|\back N(e)\right|\,\!.$$

Moreover, given a~$k$-graph~$H$ and two disjoint $(k-1)$-sets of vertices $e,f\in V(H)^{(k-1)}$,  we write~$e \triangleright f$ to mean $e\in \back N(f)$\,.
Thus, it is easy to see that~$Z_\ell$ can be viewed as a sequence of $(k-1)$-sets of vertices $e_1,\dots, e_\ell$ such that~$e_i\triangleright e_{i+1}$ for every $i\in [\ell]$ (where the sum is taken modulo $\ell$).

We split the proof in the lower and upper bound.

\subsection{Upper bound}

Here we prove the following lemma that yields the upper bound in Theorem~\ref{thm:main2}.
\begin{lemma}\label{lem:upper}
    Let~$k\geq 3$.
    For every~$d\in (0,1]$, there is a positive integer~$\ell\in \mathds N$ such that
    $$\gamma (Z_\ell) \leq d\,.$$
\end{lemma}

We will make use of the following lemma due to Mubayi and Zhao~\cite{MZ:07}.

\begin{lemma}\label{lem:MZ}
Fix~$k\geq 2$.
Given~$\eps, \alpha>0$ with $\alpha+\eps<1$, there exists an~$m_0\in \mathds N$ such that the following holds for every~$n$-vertex~$k$-graph~$H$ with~$\delta(H)\geq (\alpha + \eps)n$.
For every integer~$m$ with~$m_0\leq m\leq n$, the number of~$m$-sets $S\subseteq V(H)$ satisfying~$\delta(H[S])\geq (\alpha + \eps/2)m$ is at least $\tfrac{1}{2}\binom{n}{m}$.
\end{lemma}

For positive integers~$f, c$ and a~$k$-graph $F$ on $f$ vertices, denote the \emph{$c$-blow-up of~$F$} by~$F(c)$.
This is the $f$-partite~$k$-graph $F(c)=(V, E)$ with $V = V_1 \dot\cup \dots \dot\cup V_f$, $|V_i| = c$ for~$1\leq i \leq f$, and $E = \{v_{i_1}\cdots v_{i_k} : v_{i_j} \in V_{i_j} \text{ for every }j\in[k]  \text{ and } i_1,\dots, i_k \in E(F)\}$. 

By cyclically going around the vertices, it is easy to check that the blow-up of a zycle of length~$r$ contains zycles whose length is a multiple of~$r$.

\begin{fact}\label{fact:blowup}
    For~$k,r\geq 3$ and~$c\in \mathds N$, we have~$Z_{cr}\subseteq Z_r(c)$. 
\end{fact}

The following supersaturation result follows from a standard application of~Lemma~\ref{lem:MZ} combined with a classical result by Erd\H os \cite{E:64}.

\begin{prop}\label{prop:supersaturation}
    Let~$t,k,c\in \mathds N$ with~$k\geq 2$ and let~$\mathcal F=\{F_1, \dots F_t\}$ be a finite family of~$k$-graphs with~$\vert V(F_i)\vert=f_i$ for all $i\in [t]$.
    For every~$\eps>0$, there exists a~$\zeta>0$ such that for sufficiently large $n\in \mathds N$, the following holds.
    Every $n$-vertex~$k$-graph $H$ with~$\delta(H)\geq (\gamma(\mathcal F)+\eps)n$
    contains~$\zeta \binom{n}{f_i}$ copies of~$F_i$ for some~$i\in[t]$.
    Consequently, $H$ contains a copy of~$F_i(c)$.
\end{prop}

\begin{proof}
    Given~$t,k,c$ and~$\eps>0$, let~$m_0\in \mathds N$ be given by Lemma~\ref{lem:MZ}, and let~$C\in\mathds{N}$ with~$C^{-1}\ll c^{-1}$. 
    Let~$m\in\mathds{N}$ with~$m^{-1}\ll \eps, m_0^{-1}, C^{-1},f_i^{-1},k^{-1},t^{-1}$, and set
    $$\zeta=\frac{1}{2t\binom{m}{\max_i f_i}}\,.$$
    Now let~$n\in\mathds{N}$ be sufficiently large, i.e.,~$n^{-1}\ll \zeta$.
    Let~$H$ be given as in the statement of the lemma.
    Due to Lemma~\ref{lem:MZ}, at least~$\frac{1}{2}\binom{n}{m}$ induced $m$-vertex subhypergraphs of $H$ have minimum codegree at least~$(\gamma(\mathcal F)+\eps/2)m$.
    Since~$m$ is sufficiently large, each of those subgraphs will contain a copy of a hypergraph in~$\mathcal F$. 
    Therefore, there exists an~$i\in[t]$ such that there are at least~$\frac{1}{2t}\binom{n}{m}$ induced~$m$-vertex subgraphs of~$H$ containing a copy of $F_i$. 

    Set~$F=F_i$ and~$f=f_i$, and define an auxiliary $f$-uniform hypergraph $G_F$ by~$V(G_F)=V(H)$ and~$E(G_F)=\{S\in V(H)^{(f)}\colon F\subseteq H[S]\}$. 
    By the counting above, we have
    $$|E(G_F)| \geq \frac{1}{2t}\frac{\binom{n}{m}}{\binom{n-f}{m-f}} = \frac{1}{2t\binom{m}{f}}\binom{n}{f}\geq \zeta \binom{n}{f}\,.$$
    A result by Erd\H os \cite{E:64} implies that~$G_F$ contains a copy of $K_{f}^{(f)}(C)$. 
    Each edge of~$K_{f}^{(f)}(C)$ corresponds to (at least) one embedding of~$F$ into $H$, in one of the at most~$f!$ possible ways that~$F$ could be embedded into the $f$ vertex classes of~$K_{f}^{(f)}(C)$ (viewed as vertex sets of~$H$).
    Thus, when colouring the edges of~$K_{f}^{(f)}(C)$ accordingly, Ramsey's theorem entails that there is a~$K_{f}^{(f)}(c)\subseteq K_{f}^{(f)}(C)$ for which all embeddings of~$F$ follow the same permutation.
    This yields a copy $F(c)$ in~$H$.
\end{proof}

No we are ready to prove Lemma~\ref{lem:upper}.

\begin{proof}[Proof of Lemma~\ref{lem:upper}]
    Given~$k\geq 3$ and~$d\in (0,1)$ (since for~$d=1$ the statement is clear), take~$t=\lceil d^{-2(k-1)}\rceil+1$ and~$\ell=(2t)!$.
    We first prove the following claim.
    \begin{clm}\label{claim:evenzycles} $\gamma(Z_2, Z_4,\dots, Z_{2t}) \leq d\,.$
    \end{clm}

    \begin{claimproof}
            
    Let~$\eps\ll 1/k,1/t,1-d$ and pick~$n\in\mathds{N}$ with $n^{-1}\ll \eps$.
    Let~$H=(V,E)$ be a~$k$-graph on~$n$ vertices with~$\delta(H)\geq (d + \eps)n$. 
    We shall prove that~$Z_{2r}\subseteq H$ for some $r\in \{1,\dots, t\}$. 
    To this end, we find a sequence of $(k-1)$-sets of vertices $e_1,\dots, e_{2r}\in V^{(k-1)}$ with~$e_i\triangleright e_{i+1}$ for every~$i\in [2r]$ (where the sum is modulo $2r$). 
    First, we show that there is a sequence of pairwise disjoint $(k-1)$-sets of vertices $e_1, e_3, \dots, e_{2t-1} \in V^{(k-1)}$ such that 
        \begin{align}\label{eq:interN}
            \big\vert
            N(e_{2i-1})^{(k-1)}
            \cap
            \back{N}(e_{2i+1}) 
            \big\vert 
            > \frac{1}{t-1} \binom{n}{k-1}+t(k-1)n^{k-2}\,, 
        \end{align}
    for every~$i\in [t-1]$.
    

        Pick $e_1$ arbitrarily.
        We choose~$e_3,\dots, e_{2t-1}$ iteratively as follows. 
        Suppose that for~$j\in[t-1]$, we have already found a sequence~$e_1,\dots, e_{2j-1}$ satisfying~\eqref{eq:interN} for every~$i\leq j$.
        Let~$U_j=\bigcup_{i\in [j]} e_{2i-1}$ and note that~$|U_j|\leq (k-1)t\leq \tfrac{\eps n}{2}$.
        The following identity holds by a double counting argument, and the inequality follows from the minimum codegree condition
        \begin{align*}
        \sum_{e\in (V\setminus U)^{(k-1)}} 
        \nqquad|N(e_{2j-1})^{(k-1)} \cap \back N(e)| 
        =\nquad
        \sum_{e\in N(e_{2j-1})^{(k-1)}} \binom{|N(e)\setminus U|}{k-1} 
        \geq 
        \binom{(d+\tfrac{\eps}{2})n}{k-1}^2\,.
        \end{align*}
    Therefore, by averaging there is an~$e_{2j+1}\in (V\setminus U_j)^{(k-1)}$ such that 
     \begin{align*}
    \big| N(e_{2j-1})^{(k-1)} \cap\back N(e_{2j+1})\big|
    \geq 
    \frac{\binom{(d+\frac{\eps}{2})n}{k-1}^2}{\binom{n}{k-1}}
    \geq&\,
    \Big(d+\frac{\eps}{4}\Big)^{2(k-1)}\binom{n}{k-1} \\[1.7pt]
    \geq&\,
    d^{2(k-1)}\binom{n}{k-1}+t(k-1)n^{k-2}\\[1.7pt]
    \geq&\,
    \frac{1}{t-1} \binom{n}{{k-1}}+t(k-1)n^{k-2}\,. 
    \end{align*}
    Hence, after~$t$ steps we found~$e_1,e_3,\dots,e_{2t-1}\in V^{(k-1)}$ satisfying~\eqref{eq:interN} for every~$i\in [t-1]$.

    Note that the number of~$(k-1)$-sets containing at least one vertex in~$\bigcup_{i\in[t]}e_{2i-1}$ is at most~$t(k-1)n^{k-2}$.
    Thus, because of~\eqref{eq:interN}, the pigeonhole principle implies that there are indices~$i,j\in[t-1]$ with~$i<j$ and~$e_{2i}\in \bigcap_{s\in\{i,j\}}\Big(N(e_{2s-1})^{(k-1)}\cap \back N(e_{2s+1})\Big)$ such that~$e_{2i}$ is disjoint from each of~$e_1,e_3,\dots,e_{2t-1}$.
    In particular, we have
    \begin{align}\label{eq:specialedge}
        e_{2i}\triangleright e_{2i+1} \qqand e_{2j-1}\triangleright e_{2i}\,.
    \end{align}
    
    Next we choose the other $(k-1)$-sets with even indices in the sequence forming~$Z_{2r}$. 
    We shall choose $j-i-1$ pairwise disjoint $(k-1)$-sets~$e_{2i+2}, \dots, e_{2j-2}\in V^{(k-1)}$ such that~$e_{2m}\in N(e_{2m-1})\cap \back N(e_{2m+1})$ for every $i<m<j$ (note that if~$j=i+1$, we are done). 
    In other words, for $i<m<j$, we need
    \begin{align}\label{eq:step}
        e_{2m-1} \triangleright e_{2m} \triangleright e_{2m+1}\,.
    \end{align}
    Moreover, the~$e_{2m}$ have to be disjoint from the already chosen sets in the sequence.
    Each set~$e \in V(H)^{(k-1)}$ can intersect at most $(k-1)n^{k-2}$ other elements of~$V(H)^{(k-1)}$.
    Thus, we can greedily pick disjoint the even sets $e_{2m}\in N(e_{2m-1})^{(k-1)} \cap \back N(e_{2m+1})$ one by one for each~$i<m<j$.
    Indeed, for every~$m\leq j-i-1$, the number of~$(k-1)$-sets in $N(e_{2m-1})^{(k-1)}\cap \back N(e_{2m+1})$ which do not intersect any previously chosen~$(k-1)$-set in the sequence is at least
    $$|N(e_{2m-1})^{(k-1)}\cap \back N(e_{2m+1})| - 2t(k-1)n^{k-2} 
    \overset{\eqref{eq:interN}}{\geq}
    \frac{1}{t-1} \binom{n}{k-1}- t(k-1)n^{k-2} 
    \overset{\phantom{\eqref{eq:interN}}}{>} 0\,.$$
    This means that we can always pick an~$e_{2m}\in N(e_{2m-1})^{(k-1)} \cap \back N(e_{2m+1})$ that is disjoint from all previously chosen sets.
    
    Putting~\eqref{eq:specialedge} and~\eqref{eq:step} together yields that the $(k-1)$-sets~$e_{2i},e_{2i+1},\dots, e_{2j-1},$ form a zycle of length~$2(j-i)\leq 2t$.  
    This concludes the proof of the claim.
    \end{claimproof}

    Let~$0<\eps\ll 1/\ell$ and~$m\geq \ell/2$. 
    Let~$n\in\mathds{N}$ with~$n^{-1}\ll \eps$ and let~$H$ be a~$k$-graph with~$\delta(H)\geq (d+\eps)n$.
    We shall prove that~$Z_\ell\subseteq H$. 
    Notice that Proposition~\ref{prop:supersaturation} and Claim~\ref{claim:evenzycles} imply that~$H$ contains a copy of~$Z_{2r}(m)$ with~$r\in \{1,\dots, t\}$. 
    Applying Fact~\ref{fact:blowup} with $c = \tfrac{\ell}{2r}\leq m$, we obtain a copy of~$Z_\ell$ in $H$ as desired.
\end{proof}

\subsection{Lower bound}\label{sec:lower}

The following construction will provide an example of a~$Z_\ell$-free hypergraph with large minimum codegree.

\begin{dfn}\label{dfn:alggraph}
    Let~$n,p,k\in \mathds N$ be such that $p$ is a prime,~$k \geq 2$ and~$p\mid n$.
    We define the $n$-vertex~$k$-graph~$\mathds F_p^{(k)}(n)$ as follows.
    The vertex set consists of~$p$ disjoint sets of size~$\frac{n}{p}$ each, i.e., $V(\mathds F_p^{(k)}(n)) =  V_0 \dcup \dots \dcup V_{p-1}$ with~$\vert V_i\vert=\frac{n}{p}$ for all~$i\in [p]$.
    Given a vertex~$v\in V(\mathds F_p^{(k)}(n))$ we write~$\mathfrak f(v)=i$ if and only if~$v\in V_i$ for~$i\in \{0,1,\dots, p-1\}$.
    We define the edge set of~$\mathds F_p^{(k)}(n)$~by
    $$v_1\cdots v_k \in E(\mathds F_p^{(k)}(n)) \Leftrightarrow 
    \begin{cases}
        \mathfrak f(v_1)+\dots+ \mathfrak f(v_k) \equiv 0 \bmod p \text{ and } \mathfrak f(v_i) \neq 0 \text{ for some }i\in [k] \text{, or}\\
        \mathfrak f(v_{\sigma(1)})=\dots=\mathfrak f(v_{\sigma(k-1)})=0 \text{ and } \mathfrak f(v_{\sigma(k)})=1 \text{ for some }\sigma\in S_k\,.
    \end{cases}$$
\end{dfn}

When~$k$ is obvious from the context, we omit it from the notation and we always consider the indices of the clusters modulo $p$.

\pagebreak

\begin{lemma}\label{lem:lower}
    Let~$k\geq 3$. For every~$\ell\geq 2$, we have $\frac{1}{2(k-1)^{\ell}} \leq \gamma(Z_\ell^{(k)})$.
\end{lemma}

\begin{proof}
    Given~$k\geq 3$ and~$\ell\geq 2$, let~$n,p\in \mathds N$ be such that~$p\mid n$, $p$ is a prime larger than~$k$, and $n^{-1}\ll p^{-1}<\frac{1}{(k-1)^\ell+1}$. 
    Observe that by the Bertrand–Chebyshev theorem we might take $p\leq 2(k-1)^{\ell}$. 
    We shall prove that
    \begin{align}\label{eq:toprove}
      \delta(\mathds F_p(n)) = \tfrac{n}{p} \geq \frac{n}{2(k-1)^\ell}\qqand 
      Z_\ell \not\subseteq \mathds F_p(n)\,.
    \end{align}
    To check the codegree condition in \eqref{eq:toprove}, take a $(k-1)$-set  of vertices~$v_1,\dots,v_{k-1}$. 
    If there is an~$i\in [k-1]$ such that $\mathfrak f (v_i) \neq 0$, then let~$j$ be the only solution in $\{0,1,\dots, p-1\}$ to the equation
    $$\mathfrak f(v_1)+\dots+ \mathfrak f(v_{k-1}) + x \equiv 0 \pmod p\,.$$
    Then,~$N(v_1\cdots v_{k-1})\supseteq V_j$ and therefore~$d(v_1\cdots v_{k-1})\geq \tfrac{n}{p}$.
    If~$f(v_i)=0$ for all $i\in [k-1]$, then~$N(v_1\cdots v_{k-1})=V_1$ and we obtain~$d(v_1\cdots v_{k-1})=\tfrac{n}{p}$. 

    To check the second part of~\eqref{eq:toprove}, assume that there are~$r\geq 2$ and sets $e_1,\dots, e_r \in V(\mathds F_p(n))^{(k-1)}$ forming a copy of~$Z_r$, i.e., we have~$e_i\triangleright e_{i+1}$ for all~$i$.
    Here, and for the rest of the proof, we take the sum of indices of the~$e_i$'s to be modulo $r$.
    We shall prove that~
    \begin{align}\label{toprove1}
    r>\ell\,.     
    \end{align}
    
    The following claim states that there is an~$i_0$ for which~$e_{i_0}$ is completely contained in one of the clusters of~$\mathds F_{p}(n)$. 
    Moreover, that cluster is not~$V_0$.
    \begin{clm}\label{claim:trapped}
        There is an~$i_0\in [r]$ and a~$j\in[p-1]$ such that~$e_{i_0}\subseteq V_j$.
    \end{clm}
    \begin{claimproof}
        Fix any~$i\in [r]$, let~$e_i=v_1\cdots v_{k-1}$, and pick~$v_k \in e_{i+1}$ arbitrarily.
        We consider four cases. 
        \begin{enumerate}[label=\caselabel,     
        leftmargin=\widthof{[Case (4)]}+\labelsep, 
        labelindent=0pt]
        \item $|e_i\cap V_0|=k-1$.
        \label{case:k-1}
        \smallskip
        
        By Definition~\ref{dfn:alggraph} and since~$v_1\cdots v_k\in E(\mathds F_p(n))$, we have~$v_k\in V_1$. 
        Since we picked~$v_k\in e_{i+1}$ arbitrarily, we have that $e_{i+1} \subseteq V_1$ and finish the proof of this case by taking~$i_0=i+1$.
        
        \item $|e_i\cap V_0| < k-2$. 
        \smallskip
        
        Let~$j\equiv -(\mathfrak f(v_1)+\cdots + \mathfrak f(v_{k-1})) \bmod p$.
        By Definition~\ref{dfn:alggraph} and since~$v_1\cdots v_k\in E(\mathds F_p(n))$, we have
        $$0\equiv
        \mathfrak f(v_1)+\cdots + \mathfrak f(v_k) 
        \equiv \mathfrak f(v_k) - j\,.$$
        This means that~$v_k\in V_j$ and since we picked~$v_k\in e_{i+1}$ arbitrarily, similarly as above we get~$e_{i+1}\subseteq V_j$. 
        If~$j\not\equiv 0$, we finish by taking~$i_0=i+1$.
        If~$j\equiv0$, the claim follows from Case~\hyperref[case:k-1]{(\!\textit{1}\!)} for~$e_{i+1}$ instead of~$e_i$.

        \smallskip
        \item $|e_i\cap V_0| = k-2$ and $|e_i\cap V_1| =0$. \label{case:=k-2andno1} 
        \smallskip
        
        This case follows from similar arguments as the previous one.

        \smallskip
        \item $|e_i\cap V_0| = k-2$ and $|e_i\cap V_1| =1$.
        \smallskip 
        
        By Definition~\ref{dfn:alggraph}, we either have~$v_k\in V_0$ or~$v_k\in V_{p-1}$.
        Thus, since we picked~$v_k\in e_{i+1}$ arbitrarily, we certainly have~$e_{i+1}\subseteq V_0\cup V_{p-1}$. 
        Hence,~$|e_{i+1}\cap V_1|=0$ and so the proof follows from Cases~\hyperref[case:k-1]{(\!\textit{1}\!)} - \hyperref[case:=k-2andno1]{(\!\textit{3}\!)} for~$e_{i+1}$ instead of~$e_i$.
        \end{enumerate}
    \vspace{-10pt}
    \end{claimproof}

    We now show that for every~$i\in[r]$,
        \begin{align}\label{toprove2}
            \text{if }e_i \subseteq V_{j} \text{ with $j\not\equiv 0 \bmod p$, then }e_{i+1}\subseteq V_{(1-k)j}\,.
        \end{align}
    Indeed, let~$e_{i}=v_1\cdots v_{k-1}\subseteq V_{j}$ and pick~$v_k\in e_{i+1}$ arbitrarily. 
    Since~$\mathfrak f(v_i) \equiv j \bmod p$ for~$i\in [k-1]$, we have
    $$\mathfrak f(v_1)+\dots+ \mathfrak f(v_{k-1}) \equiv (k-1)j \pmod p\,.$$
    Therefore, since~$e_i\triangleright e_{i+1}$ implies~$v_1\cdots v_k\in E(\mathds F_p(n))$ and because~$\mathfrak f(v_i) \equiv j\not\equiv 0 \bmod p$ for~$i\in [k-1]$, we have
    $$0\equiv \mathfrak f(v_1)+\dots+ \mathfrak f(v_{k}) \equiv (k-1)j + \mathfrak f(v_k) \pmod p\,.$$
    Hence~$\mathfrak f (v_k) \equiv (1-k)j$, meaning that~$v_k\in V_{(1-k)j}$.
    Since we picked~$v_k\in e_{i+1}$ arbitrarily, we have~$e_{i+1}\subseteq V_{(1-k)j}$ proving~\eqref{toprove2}.

    Finally, we are ready to show~\eqref{toprove1}.
    Let~$i_0$ and~$j$ be given by~Claim~\ref{claim:trapped}.
    As~$p$ is a prime,~$\mathds F_p$ is a field.
    Together with~$j\nequiv 0$, this entails that~$(1-k)^sj\nequiv 0 \pmod p$ for all~$s\in[r]$.
    Thus,~$r$ applications of \eqref{toprove2} imply that 
    $$e_{i_0+r} \subseteq V_{m} \text{ with }m\equiv (1-k)^r j \pmod p\,.$$
    Since~$e_{i_0+r} = e_{i_0} \in V_j$, we have $(1-k)^rj \equiv j \pmod p$, and as~$j\nequiv 0$, we have~$(1-k)^r \equiv 1$. 
    Recalling that we chose $p$ such that $p>(k-1)^{\ell}+1$,~\eqref{toprove1} follows. 
\end{proof}

\section{Proof of Theorem~\ref{thm:min}}

\subsection{Method}

As mentioned in the introduction, to prove Theorem~\ref{thm:min} we apply the method developed by the authors together with Sales in \cite{PSS}. 

\begin{dfn}\label{dfn:picture}
Given a~$k$-graph~$H=(V,E)$, a \emph{picture} is a tuple~$(v,m,\mathcal{L},\mathcal{B})$, where
\begin{enumerate}[label=\rmlabel]
    \item $v\in V$,
    \item $m\in\mathds{N}$,
    \item $\mathcal L$ is a collection of~$m$-tuples~$\mathcal{L}\subseteq (V\setminus\{v\})^m$, and 
    \item $\mathcal{B}\subseteq [m]^{(k-1)}$ is a fixed family of~$(k-1)$-subsets of~$V(H)$, 
\end{enumerate}
such that for every~$(x_1,\dots,x_m)\in \mathcal{L}$ and every~$i_1\cdots i_{k-1}\in\mathcal{B}$, the~$k$-sets ~$vx_{i_1}\cdots x_{i_{k-1}}$ are edges of~$H$. 
That is to say, $x_{i_1}\cdots x_{i_{k-1}}$ is an edge in the link of~$H$ at~$v$.
\end{dfn}

We use pictures to find a copy of a $k$-graph $F$ on~$H$.
Roughly speaking, we say that a picture is \emph{nice} if it `encodes' a set of edges that would yield a copy of~$F$, but whose existence we cannot (yet) guarantee when considering the link of~$H$ at~$v$.

\begin{dfn}\label{dfn:nicepic}
    Given~$k$-graphs~$F$ and~$H=(V,E)$, and vertex set~$S\subseteq V$, we say that a picture~$(v,m,\mathcal{L},\mathcal{B})$ is \emph{$S$-nice for~$F$}, if for every~$w\in S$ and every~$(x_1,\dots,x_m)\in\mathcal{L}$, the hypergraph with vertex set~$V$ and edge set $$E\cup \{wx_{i_1}\cdots x_{i_{k-1}}:i_1\cdots i_{k-1}\in\mathcal{B}\}$$ contains a copy of~$F$.
\end{dfn}

If~$F$ is clear from the context, we speak simply of~$S$-nice pictures.
The following lemma describes how the existence of~$S$-nice pictures implies that~$H$ contains a copy of~$F$.

\begin{lemma}\label{lem:method}
    Let~$F$ be a $k$-graph.
    Given~$\xi,\zeta> 0$ and~$c,m\in\mathds N$, let~$n\in\mathds{N}$ such that~$n^{-1}\ll \xi, \zeta, \vert V(F)\vert^{-1},c^{-1},m^{-1}$, and let~$H$ be an~$n$-vertex~$k$-graph.
    
    Suppose that there are~$m\in\mathds{N}$ and~$\mathcal{B}\subseteq[m]^{(k-1)}$ such that for every~$S\subseteq V(H)$ with~$\vert S\vert\geq c$, there is an~$S'$-nice picture~$(v,m,\mathcal L,\mathcal{B})$, with $v\in S$, $S'\subseteq S$,~$|S'|\geq \xi |S|$, and~$|\mathcal L|\geq \zeta n^m$. 
    Then $H$ contains a copy of~$F$. 
\end{lemma} 

\begin{proof}
    Let~$t=\lceil\zeta^{-1}\rceil+1$. 
    By iteratively applying the conditions of the lemma, we find a nested sequence of subsets~$V(H)=S_0\supseteq S_1\supseteq \dots \supseteq S_t$ such that for~$i\in[t]$, there are $S_{i}$-nice pictures $(v_i,m, \mathcal L_i,\mathcal{B})$ satisfying~$v_i\in S_{i-1}$, $|S_{i}|\geq \xi^{i} n>c$, and $|\mathcal L_i| \geq \zeta n^m$.

    Since~$t\geq \zeta^{-1}+1$, by the pigeonhole principle, there are two indices~$0<i<j\leq t$ such that~$\mathcal L_i\cap \mathcal L_j\neq \emptyset$. 
    Let~$(x_1,\dots,x_m)\in \mathcal L_i\cap \mathcal L_j$. 
    Then because $(v_i,m,\mathcal L_i,\mathcal{B})$ is an $S_i$-nice picture and~$v_j\in S_{j-1}\subseteq S_i$, Definition~\ref{dfn:nicepic} guarantees that
    $$E(H)\cup \{v_jx_{i_1}\cdots x_{i_{k-1}} \colon i_1\cdots i_{k-1}\in \mathcal{B}\}$$ 
    contains a copy of~$F$.
    Since~$(v_j,m,\mathcal L_j,\mathcal{B})$ is a picture, Definition~\ref{dfn:picture} yields~$v_jx_{i_1}\cdots x_{i_{k-1}}\in E(H)$ for all~$i_1\cdots i_{k-1}\in\mathcal{B}$.
    Thus, we conclude that this copy of~$F$ is in fact in~$H$.
\end{proof}

Now we apply Lemma~\ref{lem:method} to prove Theorem~\ref{thm:min}.

\subsection{Proof of Theorem~\ref{thm:min}}
    Let~$\ell\geq 3$ be an integer and let~$\eps>0$.
    Let~$\xi,\zeta > 0$, and let~$n, c\in \mathds N$ such that~$n^{-1}\ll c^{-1}\ll \zeta,\xi\ll\eps$. 
    Let~$H$ be a~$3$-graph with $\delta(H)\geq \eps n$.
    We aim to show that~$Z_{\ell}^-\subseteq H$.
    Set~$m=2$ and~$\mathcal{B}=\big\{\{1,2\}\big\}$, then due to Lemma~\ref{lem:method}, we only need to prove that for every~$S\subseteq V(H)$ of size at least~$c$, there is an $S'$-nice picture~$(v,2, \mathcal L,\big\{\{1,2\}\big\})$ with~$v\in S$,~$S'\subseteq S$, $|S'|\geq \xi |S|$, and~$|\mathcal L|\geq \zeta n^{2}$.

    Given~$S\subseteq V(H)$ with~$\vert S\vert\geq c$, take any vertex~$v\in S$ and let~$V=V(H)\setminus\{v\}$.
    Observe that using the minimum codegree condition and the above hierarchy, we have
    \begin{align}\label{eq:doublecounting}
        \sum_{bb'\in V^{(2)}} |N_{L_v}(b)\cap N_{L_v}(b')\cap S| = \sum_{u\in S\setminus \{v\}} \binom{d_{L_v}(u)}{2} \geq \binom{\eps n}{2}(|S|-1) \geq \xi \binom{n}{2}|S|\,,
    \end{align}
    where~$L_v$ denotes the link of~$H$ at~$v$.
    Thus, by averaging there is a pair~$b_1,b_2\in V$ such that $\vert N_{L_v}(b_1)\cap N_{L_v}(b_2)\cap S\vert\geq\xi\vert S\vert$.
    We pick~$S'\subseteq N_{L_v}(b_1)\cap N_{L_v}(b_2)\cap S$ with~$\vert S'\vert=\lceil\xi\vert S\vert\rceil$.

    Since~$\delta(H) - 2\ell - |S'|\geq \eps n/2\geq 2$ we can greedily pick pairwise disjoint pairs of vertices~$e_1,\dots, e_{\ell-2}\in (V(H)\setminus S')^{(2)}$ such that
    \begin{align}\label{eq:greedy}
        b_1b_2=e_1\triangleright e_2 \triangleright \cdots \triangleright e_{\ell-2}\,.
    \end{align}
    Now let~$R=\bigcup_{i\in [\ell-2]}e_i$ and take
    $$\mathcal L = \{(x_1,x_2) \in V^2\colon x_1\in N_H(e_{\ell-2})\setminus R \text{ and } x_2\in N_H(x_1v)\setminus R\}\,.$$
    Note that~$|\mathcal L| \geq (\delta(H)/2)^2\geq \eps^2 n^2/4\geq\zeta n^2$.
    Further, since~$x_1x_2\in E(L_v)$ for every~$(x_1,x_2)\in\mathcal{L}$,~$(v,m,\mathcal L,\mathcal{B})$ is a picture in~$H$.
    Moreover, observe that it is~$S'$-nice. 
    Indeed, we only we need to check that for any~$u\in S'$ and~$(x_1,x_2)\in\mathcal L$, the hypergraph with edges~$E(H)\cup \{ux_1x_2\}$ contains a copy of~$Z_{\ell}^-$.
    For this, note that in~$E(H)\cup \{ux_1x_2\}$ we have $x_1x_2\triangleright uv$.
    Further,~$u\in S'$ and the choice of~$b_1$ and~$b_2$ imply~$uv\triangleright b_1b_2$.
    Together with~\eqref{eq:greedy}, this gives~$x_1x_2\triangleright uv\triangleright e_1\triangleright \cdots \triangleright e_{\ell-2}$, and using the fact that~$x_1\in N(e_{\ell-2})$, we obtain a copy of~$Z_{\ell}^-$ (where the missing edge is~$x_2e_{\ell-2}$).

\section{Concluding Remarks}

Following a very similar proof as that for Theorem~\ref{thm:min}, we can show a general upper bound for~$\gamma(Z_\ell^{(3)})$ for every~$\ell\geq 3$.
\begin{prop}\label{prop:1/2} For~$\ell\geq 3$,
    $\gamma(Z_\ell^{(3)})\leq 1/2$.
\end{prop}

\begin{proof}
    Given~$\ell\geq 3$ and~$\eps>0$, let~$\xi,\zeta > 0$ and~$n, c\in \mathds N$ such that~$n^{-1}\ll c^{-1}\ll \zeta,\xi\ll\eps$. 
    Let~$H$ be a~$3$-graph with $\delta(H)\geq \big(\frac{1}{2}+\eps\big) n$.
    We aim to show that~$Z_{\ell}\subseteq H$.
    As in the proof of Theorem~\ref{thm:min}, we pick $m=2$ and $\mathcal B=\big\{\{1,2\}\big\}$ and due to Lemma~\ref{lem:method}, we only need to prove that for every~$S\subseteq V(H)$ of size at least~$c$, there is an~$S'$-nice picture~$(v,2, \mathcal L,\big\{\{1,2\}\big\})$ with~$v\in S$,~$S'\subseteq S$, $|S'|\geq \xi |S|$, and~$|\mathcal L|\geq \zeta n^{2}$. 
    
    For the first part of the proof we proceed as in the proof of Theorem~\ref{thm:min} and we only use~$\delta(H)\geq \eps n$. 
    In particular, we obtain two vertices~$b_1,b_2\in V(H)\setminus \{v\}=:V$ and a set~$S'\subseteq N_{L_v}(b_1)\cap N_{L_v}(b_1)\cap S$ with $|S'|=\lceil \xi |S|\rceil$.
    Moreover, we again greedily pick pairwise disjoint pairs of vertices~$e_1,\dots, e_{\ell-2}\in (V\setminus S')^{(2)}$ satisfying~\eqref{eq:greedy}.
    The set~$\mathcal L$ is chosen differently.
    Set~$R=\bigcup_{i\in[\ell-2]}e_i$ and 
    \begin{align}\label{eq:defL}
        \mathcal L = \{(x_1,x_2)\in V^2\colon x_1,x_2\in N(e_{\ell-2})\setminus R \text{ and }x_1x_2\in E(L_v)\}\,.
    \end{align}
    Observe that given~$x_1\in N(e_{\ell-2})\setminus R$, any vertex~$x_2\in (N(xv)\cap N(e_{\ell-2}))\setminus R$, gives rise to~$(x_1,x_2)\in \mathcal L$. 
    Furthermore, since~$\delta(H)\geq (1/2+\eps)n$,
    $$|(N(xv)\cap N(e_{\ell-2}))\setminus R| \geq \eps n - 2\ell \geq \frac{\eps}{2}n\,,$$
    and similarly we have~$N(e_{\ell-2})\setminus R\geq n/2$.
    Therefore, we obtain~$|\mathcal L|\geq \eps n^2/4$, and since~$x_1x_2\in E(L_v)$ for all~$(x_1,x_2)\in\mathcal{L}$,~$(v,m,\mathcal{L},\mathcal{B})$ is a picture in~$H$.

    To see that the tuple~$(v,m,\mathcal{L},\mathcal{B})$ is indeed an $S'$-nice picture, we shall prove that for every~$u\in S'$ and~$(x_1,x_2)\in \mathcal L$, the hypergraph with (vertex set~$V(H)$ and) edges~$E(H)\cup \{ux_1x_2\}$ contains a copy of~$Z_{\ell}$.
    Indeed, the definition of~$\mathcal{L}$ implies $x_1x_2v\in E(H)$ and therefore~$x_1x_2\triangleright uv$ in $E(H)\cup \{ux_1x_2\}$.
    Also due to the definition of~$\mathcal{L}$, we have~$x,y\in N(e_{\ell-2})$ and thus,~$e_{\ell-2}\triangleright xy$.
    Moreover,~$u\in S'$ and the choice of~$b_1$ and~$b_2$ entails~$uv\triangleright b_1b_2=e_1$.
    Combining this with~\eqref{eq:greedy}, we obtain~$uv\triangleright e_1\triangleright\dots e_{\ell-2}\triangleright x_1x_2\triangleright uv$, that is a copy of~$Z_\ell$, in~$E(H)\cup \{ux_1x_2\}$.
\end{proof}

It would be interesting to know whether Proposition~\ref{prop:1/2} is sharp for some~$\ell\geq 3$.
The following construction gives a lower bound of~$1/3$ for the codegree Tur\'an density of any zycle of length not divisible by~$3$.
Let~$n\in\mathds N$ be divisible by~$3$ and let~$H=(V,E)$, where~$V=V_1\dcup V_2 \dcup V_3$ with~$|V_i|=n/3$ and~$E=\{uvw\in V^{(3)}\colon u,v\in V_i \text{ and }w\in V_{i+1}\}$, where the sum is taken modulo~$3$. 
It is not hard to check that~$\delta(H)\geq n/3$ and that~$Z_\ell\not\subseteq H$ for every~$\ell$ not divisible by~$3$.

Observe that~$Z_2^{(3)}=K_4^{(3)}$.
For this $3$-graph, a well-known conjecture by Czygrinow and Nagle~\cite{CN:01} states that~$\gamma(Z_2^{(3)}) = \gamma(K_4^{(3)}) = 1/2$. 
Regarding the next case,~$Z_3^{(3)}$, note that its codegree Tur\'an density is not bounded by the previous construction. 
The following~$3$-graph entails~$\gamma(Z_3^{(3)})\geq 1/4$, and in fact it provides the same lower bound for every~$Z_\ell^{(3)}$ with~$\ell$ not divisible by~$4$.
Let~$n\in \mathds N$ divisible by~$4$ and let~$H=(V,E)$, where $V=V_1\dcup V_2\dcup V_3\dcup V_4$ with~$|V_i|=n/4$.
Define the edges of~$H$ as 
$$E=\{xyz \colon \text{$x, y\in V_i$ and~$z\in V_{i+1}$}\} \dcup \{xyz\colon x\in V_1, y\in V_2, z\in V_3\cup V_4\}\,,$$
where the sum of indices is taken modulo~$4$.
Clearly,~$\delta(H)\geq n/4$. 
To see that~$Z_\ell\not\subseteq H$ for~$\ell$ not divisible by~$4$, it can be checked that all zycles are of the form $e_1\triangleright \dots \triangleright e_r$ such that~$e_i\subseteq V_{j_i}$ for some~$j_i\in[4]$. 
Together with Proposition~\ref{prop:1/2}, this yields $$\frac{1}{4}\leq \gamma(Z_3^{(3)})\leq \frac{1}{2}\,.$$

\begin{prob}
    Determine the value of~$\gamma(Z_3^{(3)})$.
\end{prob}

\bigskip

On a different note, recall that Theorem~\ref{thm:min} states that~$Z_\ell^{(3)}$ is (inclusion) minimal with respect to the property of having strictly positive codegree Tur\'an density. 
It would be interesting to know if this also holds for larger uniformities. 

\begin{quest}
    For~$k>3$ and sufficiently large~$\ell$, what are the minimal subgraphs~$F\subseteq Z_\ell^{(k)}$ with~$\gamma(F)>0\,$?
\end{quest}
Let~$e_1,\dots,e_\ell$ be pairwise disjoint~$(k-1)$-sets of vertices. 
Consider the~$k$-graph whose edges are given by~$e_1\triangleright \dots \triangleright e_\ell$ plus one additional edge of the form~$e_\ell \cup\{v\}$ with~$v\in e_1$. 
Following the same arguments as in the proof of Theorem~\ref{thm:min}, we obtain that this~$k$-graph has codegree Tur\'an density zero for every~$\ell\geq 3$. 
  
\bigskip

In order to prove that the lower bound of Lemma~\ref{lem:lower} in Subsection~\ref{sec:lower}, we introduce the~$k$-graphs $\mathds F_p^{(k)}(n)$ that have large minimum codegree and are~$Z_\ell^{(k)}$-free for small~$\ell$. 
It would be interesting to study the codegree Tur\'an density of $\mathds F_p^{(k)}(n)$ itself. 
Observe however, that for~$n\geq pk$ we have~$K_{k+1}^{(k)-}\subseteq \mathds F_p^{(k)}(n)$, which suggests that this problem might be very difficult for general~$n$. 

It is perhaps more natural to study the codegree Tur\'an density of the following~$k$-graph.
For~$p> k$, let $\widetilde{\mathds F}_p^{(k)}$ be the~$k$-graph on~$p(k-1)$ vertices with~$V(\widetilde{\mathds F}_p^{(k)}) = V_1\dcup \dots \dcup V_p$ where~$|V_i|=k-1$ for every~$i\in [p]$ and whose edges are given by
$$v_1\cdots v_k \in E(\widetilde{\mathds F}_p^{(k)}) 
\, \Longleftrightarrow \, 
\mathfrak f(v_1)+\dots+ \mathfrak f(v_k) \equiv 0 \bmod p\,,$$
where the function~$\mathfrak f\colon V(\widetilde{\mathds F}_p^{(k)})\to [p]$ is analogous as in Definition~\ref{dfn:alggraph}. 

\begin{prob}
    For~$k\geq 3$, determine the codegree Tur\'an density of $\widetilde{\mathds F}_p^{(k)}$. 
\end{prob}

Consider the indices of the clusters~$V_1, \dots, V_p$ of~$\widetilde{\mathds F}_p^{(k)}$ to be modulo $p$.
Observe for~$j\in [p]$, we have~$V_j\cup\{v\}\in E(\widetilde{\mathds F}_p^{(k)})$ for every~$v\in V_{(1-k)j}$. 
It follows that 
$$V_1 \triangleright V_{1-k} \triangleright \dots \triangleright V_{(1-k)^{p-2}}\triangleright V_{(1-k)^{p-1}} = V_1\,,$$
where the last identity is given by Fermat's little theorem. 
Hence, there is an~$\ell\leq p-1$ such that $Z_{\ell}\subseteq \widetilde{\mathds F}_p^{(k)}$ and therefore, Lemma~\ref{lem:lower} yields~$\gamma(\widetilde{\mathds F}_p^{(k)})\geq \tfrac{1}{2(k-1)^p} > 0$. 

\begin{quest}
    For~$k\geq 3$, is it true that~$\displaystyle\lim_{p\to \infty} \gamma(\widetilde{\mathds F}_p^{(k)})=0\,$?
\end{quest}

\printbibliography

\end{document}